\documentclass[oneside,leqno,9pt]{amsart}

\newcommand{\mytitle}{\LARGE\normalfont\scshape Moduli Problems in Abelian Categories and the Reconstruction Theorem}
\newcommand{\myname}{\normalfont\footnotesize John Calabrese \& Michael Groechenig}

\usepackage{thmtools}

\usepackage{eulervm}
\usepackage{amssymb}
\usepackage{fullpage}
\usepackage{microtype}
\usepackage[utf8]{inputenc}
\usepackage[british]{babel}
\usepackage{hyperref}
	\hypersetup{colorlinks=false,pdftitle=\mytitle,pdfauthor=\myname}
\usepackage{tikz}
	\usetikzlibrary{arrows,matrix,shapes,decorations}
\usepackage[colorinlistoftodos, backgroundcolor=blue!20, linecolor=blue!20, textwidth=1.7cm, textsize=scriptsize]{todonotes}
\usepackage{tikz-cd}

\numberwithin{equation}{section}

\declaretheoremstyle[
spaceabove=6pt, spacebelow=6pt,
headfont=\normalfont\scshape,
notefont=\normalfont,
notebraces={(}{)},
bodyfont=\normalfont,
postheadspace=.5em,
headpunct={ --},
headformat=swapnumber
]{tstyle}
\declaretheorem[numberwithin=section,style=tstyle,name=Theorem]{thm}
\declaretheorem[numbered=no,style=tstyle,name=Theorem]{thm*}
\declaretheorem[sibling=thm,style=tstyle,name=Proposition]{prop}
\declaretheorem[sibling=thm,style=tstyle,name=Lemma]{lem}
\declaretheorem[sibling=thm,style=tstyle,name=Definition]{defn}
\declaretheorem[sibling=thm,style=tstyle,name=Corollary]{cor}
\declaretheorem[numbered=no,style=tstyle,name=Corollary]{cor*}

\declaretheoremstyle[
spaceabove=6pt, spacebelow=6pt,
headfont=\small\normalfont\itshape,
notefont=\small\normalfont\itshap,
notebraces={(}{)},
bodyfont=\small\normalfont,
postheadspace=1em,
headpunct=:,
qed=$\blacksquare$
]{pstyle}
\declaretheorem[numbered=no,style=pstyle,name=Proof]{prf}

\declaretheoremstyle[
	spaceabove=6pt, spacebelow=6pt,
	headfont=\normalfont\itshape,
	notebraces={(}{)},
	bodyfont=\normalfont,
	postheadspace=.5em,
	headformat=swapnumber
]{rstyle}

\declaretheorem[sibling=thm,name=Remark,style=rstyle]{rmk}
\declaretheorem[numbered=no,name=Remark,style=rstyle]{rmk*}

	\newcommand{\cat}[1]{\mathtt{#1}}

	\newcommand{\onto}{\twoheadrightarrow}
	
	\renewcommand{\O}{\mathcal{O}}
	\newcommand{\m}{\mathfrak{m}}
	
	\newcommand{\QC}{\cat{QC}}

	\newcommand{\Set}{\cat{Set}}

\DeclareMathOperator{\Hom}{Hom}
\DeclareMathOperator{\lHom}{\underline{Hom}}
\DeclareMathOperator{\id}{id}

\DeclareMathOperator{\supp}{supp}

\usepackage{dsfont}
\usepackage{mathrsfs}
\newcommand{\qc}{\cat{qc}}
\newcommand{\Z}{\mathds{Z}}
\newcommand{\C}{\mathds{C}}

\newcommand{\G}{\mathds{G}}

\newcommand{\ltensor}{\stackrel{\text{\tiny L}}{\otimes}}
\DeclareMathOperator{\RlHom}{R\underline{Hom}}
\DeclareMathOperator{\RHom}{RHom}
\DeclareMathOperator{\Spec}{Spec}
\DeclareMathOperator{\End}{End}
\DeclareMathOperator{\Mod}{Mod}
\DeclareMathOperator{\Pt}{Pt}
\DeclareMathOperator{\Ann}{Ann}
\DeclareMathOperator{\Aut}{Aut}
\DeclareMathOperator{\Pic}{Pic}

\DeclareMathOperator{\et}{{\acute{e}t}}
\DeclareMathOperator{\Ind}{Ind}
\DeclareMathOperator{\coh}{coh}
\DeclareMathOperator{\Equiv}{Equiv}
\DeclareMathOperator{\Iso}{Iso}
\DeclareMathOperator{\Br}{Br}
\newcommand{\opp}{\text{op}}
\newcommand{\modo}[1]{{#1}\text{-}\!\Mod}
\newcommand{\Alg}{\cat{Alg}}
\newcommand{\Grpd}{\cat{Grpd}}
\newcommand{\BGm}{\text{B}\mathds{G}_\text{m}}

\newcommand{\Gm}{\mathds{G}_\text{m}}

\newcommand{\longto}{\longrightarrow}

\makeatletter
\providecommand\@dotsep{5}
\renewcommand{\listoftodos}[1][\@todonotes@todolistname]{%
  \@starttoc{tdo}{#1}}
\makeatother

\begin{document}
	\author{\myname}
	\title{\mytitle}
	\address{Imperial College London}
	\email{john.calabrese@imperial.ac.uk \\ m.groechenig@imperial.ac.uk}
	\begin{abstract}
		We give a moduli-theoretic proof of the classical theorem of Gabriel, stating that a scheme can be reconstructed from the abelian category of quasi-coherent sheaves over it.
		The methods employed are elementary and allow us to extend the theorem to (quasi-compact and separated) algebraic spaces.
		Using more advanced technology (and assuming flatness) we also give a proof of the folklore result that the group of autoequivalences of the category of quasi-coherent sheaves consists of automorphisms of the underlying space and twists by line bundles.
		We apply our strategy to prove analogous statements for categories of sheaves twisted by a $\Gm$-gerbe.
		Our methods allow us to treat even gerbes not coming from a Brauer class.
		As a pleasant consequence, we deduce a Morita theory for sheaves of abelian categories.
	\end{abstract}
	\maketitle
	\setcounter{tocdepth}{1}
	\tableofcontents
	\section*{Introduction}
In \cite{gabriel} Gabriel showed that a scheme $X$ is completely determined by its abelian category of quasi-coherent sheaves $\qc(X)$.
This implies that, for schemes $X$ and $Y$, if there is an equivalence of abelian categories $\qc(X) \simeq \qc(Y)$ then $X$ and $Y$ are isomorphic.
Gabriel's original formulation was for noetherian schemes and was later extended to arbitrary quasi-separated schemes by Rosenberg \cite[Prop 10.7.1]{rosenberg}.

In trying to generalise this theorem away from schemes, one immediately realises that it fails miserably for algebraic stacks.
For example, one may take $G = \Z/2\Z$ and its classifying stack $\text{B}G$ over $\C$.
A quick inspection shows that $\qc(\text{B}G)$ is equivalent to $\qc(X)$, where $X = \Spec \C \amalg \Spec \C$ is the disjoint union of two points.
In fact, both categories are equivalent to that of $\Z/2\Z$-graded complex vector spaces.
Clearly this breaks Gabriel's theorem, as the two are non-isomorphic.\footnote{
We must mention that if one is willing to consider $\qc(X)$ as a \emph{monoidal} category then the reconstruction theorem can be extended considerably \cite{balmer,lurie,support, bburg}.}
However, sitting in between schemes and stacks are algebraic spaces, and a natural question to ask is whether Gabriel's theorem still holds in this context.
Indeed the answer is yes, at least in the separated case.
The goal of this paper is to provide a short and easy proof of the reconstruction theorem which also applies to algebraic spaces.
\begin{thm}[continues=maincomparison]
    Let $X$ and $Y$ be quasi-compact and separated algebraic spaces over a ring $R$.
    If we have an equivalence $\qc(X) \simeq \qc(Y)$ of $R$-linear categories then there exists an isomorphism $X \simeq Y$ of algebraic spaces over $\Spec R$.\footnote{Of course in the noetherian setting one can freely interchange $\qc(X)$ with the category of coherent sheaves $\coh(X)$. This is possible as $\qc(X) = \Ind \coh(X)$, where $\Ind$ stands for the category of ind-objects \cite[Lemma 3.9]{lurie}, and $\coh(X)$ is recovered as the subcategory of \emph{compact} objects of $\qc(X)$. The trivial corollary to the theorem is then that, for noetherian and separated algebraic spaces over $R$, $\coh(X)\simeq \coh(Y) \iff X \simeq Y$.}
\end{thm}
The original proof by Gabriel was topological in spirit and cannot be ported naively to algebraic spaces.
The first step in his construction was to extract a sort of \emph{spectrum} out of an abelian category, so that from $\qc(X)$ one would recover the underlying Zariski topological space $|X|.$
This is possible as there is a correspondence between irreducible closed subsets of $X$ and certain subcategories of $\qc(X)$.

If $X$ is an algebraic space and not a scheme, however, the ringed space $(|X|,\O_X)$ is not enough to recover $X$.
A novel approach is required as one must reconstruct the functor defining $X$.
The key idea is simple: $X$ can be viewed as a subspace of \emph{points} (or \emph{pointlike} objects) of $\qc(X)$.
As we wish to recreate the functor of points of $X$ we need to know all maps $S \to X$.
In other words we need to know what \emph{families of points} of $\qc(X)$ are.

There is a natural procedure to promote $\qc(X)$ to a \emph{sheaf} or, rather, stack of abelian categories.
This is accomplished by assigning to any $S$ the category $\qc(S \times X)$.
What is crucial for our purposes is that one can make sense of this purely in terms of the category theory of $\qc(X)$, see Remark \ref{bchange}.
Within this sheaf one singles out a subsheaf of suitably defined \emph{pointlike objects} $\Pt_{\qc(X)} \subset \qc(X)$.
The main idea is that a family of pointlike objects $P \in \Pt_{\qc(X)}(S)$ over $S$ should look like the structure sheaf of the graph of a morphism $S \to X$.
Roughly, one defines a family of points over $S$ to be a quasi-coherent module $P \in \qc(S \times X)$ which is fibrewise over $S$ a skyscraper.
Categorically one phrases this by imposing a \emph{Schur}-like condition, e.g.~over a field one requires $P$ to have only trivial subobjects.
By establishing a correspondence between points in the categorical sense and graphs we are able to reassemble the functor of points of the original algebraic space.

\medskip
A more sophisticated manner of stating the fact that $\Pt_{\qc(X)}$ recovers $X$ is that $\Pt_{\qc(X)}$ is in fact $X \times \BGm$, the trivial $\Gm$-gerbe over $X.$
If one starts instead with a category of \emph{twisted sheaves} $\qc(X,\alpha)$, where $\alpha$ is a $\G_m$-gerbe, one analogously concludes that $\Pt_{\qc(X,\alpha)}$ is indeed the $\Gm$-gerbe on $X$ corresponding to $\alpha.$
Thus the general version of our main theorem provides a proof of C\u{a}ld\u{a}raru's conjecture \cite[Conjecture 4.1]{raru}, extended to algebraic spaces and arbitrary $\Gm$-gerbes.
\begin{thm}[continues=twisted]
	Let $R$ be a base ring and let $X$ and $Y$ be two quasi-compact and separated algebraic spaces over $R.$
	Let $\alpha, \beta$ be two $\Gm$-gerbes on $X$ and $Y$ respectively.
	Then if $\qc(X,\alpha) \simeq \qc(Y,\beta)$ as $R$-linear abelian categories then there exists an isomorphism $f$ of $R$-spaces between $X$ and $Y$, such that $f^*\beta = \alpha$.
\end{thm}
In \cite[Theorem 1]{perego} the theorem is proved for schemes smooth and separated over a field and Brauer classes, but it isn't shown whether the induced isomorphism $f$ carries over the gerbes.
This issue was remedied for smooth and projective varieties in \cite[Corollary 5.3]{canstel}.
Another proof of the theorem appeared in \cite{ant} which works for arbitrary quasi-compact and quasi-separated schemes and Brauer classes.
The present paper extends the range of the theorem to (separated) algebraic spaces and arbitrary $\Gm$-gerbes, and uses only elementary and underived technology.

Antithetically, by using the existence of integral kernels (which is very much a derived result), we can describe the group of autoequivalences of $\qc(X)$.
\begin{thm}[continues=auto]
	Let $X$ be a quasi-compact and separated algebraic space, \emph{flat} over a ring $R$.
	Then $\Aut_R(\qc(X)) \simeq  \Aut_R(X) \ltimes \Pic (X)$, where $\Aut_R(\qc(X))$ is the group of isomorphism classes of autoequivalences of $\qc(X)$ as an $R$-linear category and $\Aut_R(X)$ is the group of automorphisms of $X$ as a space over $\Spec R$.
\end{thm}
At least under the assumption of $X$ being a smooth and projective variety, this theorem is well-known (see for example \cite[Corollary 5.24]{huybrechts}).
This is the only result in this paper relying on any advanced machinery.
The flatness assumption is there to ensure that the product $X \times X$ needn't be derived.
We also prove an analogous result for categories of twisted sheaves on a scheme, Theorem \ref{tauto}.
In the recent preprint \cite{burgabber} the group $\Aut(\qc(X))$ is described for an arbitrary quasi-separated scheme, without any flatness assumptions and using underived methods.

Finally, we classify sheaves of abelian categories on an algebraic stack which are smooth-locally equivalent to $\qc$ (Corollary \ref{morita}).
We call these sheaves of abelian categories \emph{invertible}, by analogy with invertible sheaves of modules.
This result can be regarded as a \emph{Morita theory} for sheaves of abelian categories and fulfils the request in Remark 7.2 of \cite{ant}.
\subsubsection*{Structure of the Paper} The first two sections are \emph{stack-free} and are devoted to proving Gabriel's theorem for algebraic spaces.
In the third section we generalise the result to gerbes.
In the fourth section we use the existence of Fourier-Mukai kernels to describe the group of autoequivalences of the category of quasi-coherent sheaves.
The last section contains a minor corollary, a Morita theory for sheaves of abelian categories and a marginal remark about derived schemes.

\subsubsection*{Acknowledgements} 
\label{ssub:acknowledgements}
The first author would like to thank Tom Bridgeland for useful conversations.
The second author thanks Theo B\"uhler for pointing out a reference.
Both authors would like to warmly thank Benjamin Antieau and David Rydh for thoroughly reading an earlier manuscript and providing many helpful suggestions.
We would also like to thank Martin Brandenburg.
Finally we thank Richard Thomas for supporting a visit of the second author to London, where part of this research was developed.

The main idea for this paper was born during reading seminar on algebraic stacks held at Oxford in spring 2012.
We would like to thank all the participants.

\subsubsection*{Conventions} All rings and algebras will be commutative and unital.
Given a base ring $R$, we view an algebraic space $X$ as a functor
\begin{align*}
	X\colon R\text{-}\Alg \longto \Set
\end{align*}
with domain the category of $R$-algebras and target the category of sets (and satisfying the axioms of being an algebraic space).
Similarly, we view algebraic stacks as living inside the category of (weak) functors from $R\text{-}\Alg$ to the 2-category $\Grpd$ of groupoids.
For us a \emph{prestack} over $\Spec R$ will be any groupoid-valued functor on $R\text{-}\Alg$, i.e. not necessarily satisfying any kind of descent.

On an arbitrary site (usually the \'etale site of a space $X$) one also has prestacks and stacks of abelian categories.
However we shall informally refer to the latter also as \emph{sheaves} of abelian categories, to psychologically distinguish them from \emph{algebraic} stacks.
%
%

Given an algebraic space $X$ we will denote by $\qc(X)$ its category of quasi-coherent sheaves and by $|X|$ the underlying Zariski topological space.
Given an abelian category $\cat{C}$ we denote by $D(\cat{C})$ its unbounded derived category.
Given a ring $R$ we will write $D(R)$ for $D(R\text{-}\!\Mod)$ and given an algebraic space $X$ we will write $D(X)$ for $D(\qc(X))$.
Finally, to avoid any ambiguity with other standard definitions, we should point out that for $X$ quasi-compact with affine diagonal $D(X)$ coincides with $D_{\qc}(\O_{X_{\text{{\'e}t}}}\text{-}\Mod)$, the triangulated category of complexes of \'etale sheaves of $\O_X$-modules with quasi-coherent cohomology (modulo quasi-isomorphisms) \cite[\href{http://stacks.math.columbia.edu/tag/08h1}{Tag 08H1}]{stacks-project}.

	\section{Abelian Categories}
\subsection{Base Change}
We start by recalling the constructions needed from \cite{gaitsgory}.
Given a category $\cat{C}$ we define its \emph{centre} to be $Z(\cat{C}) = \End(1_\cat{C})$, the monoid of endomorphisms of the identity functor.
When $\cat{C}$ is additive then $Z(\cat{C})$ has the structure of a ring.
In particular, when $\cat{C} = $$R$-$\Mod$ is the category of modules over a ring $R$, then $Z(\cat{C})$ is precisely the centre of $R$.
As all the rings we care about are commutative, we will always have $Z(\modo{R}) = R$.

Let us now fix a ground ring $R$ for the remainder of this section.
Given an abelian category $\cat{C}$, an \emph{$R$-linear structure} on $\cat{C}$ consists of a morphism $R \to Z(\cat{C})$.
Unwrapping the definition shows that this is the same as a functorial action $R \to \Hom_\cat{C}(M,N)$ for any two objects $M,N \in \cat{C}$.
The prototypical example of such a structure is the following.
Let $\pi\colon X \to \Spec R$ be a scheme (or algebraic space) over $R$ and take $\cat{C} = \qc(X)$ to be the category of quasi-coherent sheaves on $X$.
Then all hom-spaces are naturally modules over $R$.

Moreover, $\qc(X)$ is also endowed with an \emph{action} of $R$-$\Mod$, in the sense that quasi-coherent sheaves on $X$ can be tensored with modules over $R$ using the pullback $\pi^*.$
In general we have the following fact.
\begin{prop}
	If $\cat{C}$ is an $R$-linear and cocomplete abelian category, there is a bifunctor (which we call \emph{pullback} or \emph{action})
	\begin{align*}
		\cat{C} \times \modo{R} \ni (E , M) &\longmapsto E \otimes_R M \in \cat{C}
	\end{align*}
	defined as follows:
	\begin{itemize}
		\item $E \otimes R  = E$,
		\item $E \otimes R^{\oplus I} = E^{\oplus I}$, for any indexing set $I$,
		\item as any module $M$ can be written as a cokernel of a morphism $R^{\oplus J} \to R^{\oplus I}$, define $E \otimes M$ to be the cokernel of $E^{\oplus J} \to E^{\oplus I}$.
	\end{itemize}
\end{prop}
\noindent The last definition is independent of the chosen presentation of $M$.
In other words $E\otimes_R (-)$ is the left adjoint of $\Hom_\cat{C}(E,-).$
We define an object $E \in \cat{C}$ to be $R$-\emph{flat} if the action functor $M \mapsto E \otimes_R M$ is exact.

Given an $R$-linear category $\cat{C}$ and a ring homomorphism $R \to R'$, we can form the \emph{base change} category $\cat{C}' = \cat{C} \otimes_R R'$.
One way to define this category is via a universal property: it is the initial cocomplete $R'$-linear abelian category admitting a colimit-preserving $R$-linear functor $\cat{C} \to \cat{C}'$.
Another way to define it, or to show that such a category exists, is as follows.
The objects of $\cat{C}'$ are given by pairs $(E,\alpha)$, 
where $E \in \cat{C}$ and
\begin{align*}
	\alpha\colon R' \otimes E \to E
\end{align*}
is such that the two natural morphisms $R' \otimes R' \otimes E \to E$ are equal.
With this model it's easy to see that there are pullback $\cat{C} \to \cat{C}'$ and forgetful $\cat{C}' \to \cat{C}$ functors, which are respectively left and right adjoints of one another.
\begin{rmk}\label{bchange}
	When $\cat{C} = \qc(X)$ is the category of quasi-coherent sheaves of an algebraic space $X$ over $R$, then the category $\qc(X) \otimes_R R'$ we abstractly defined is simply $\qc(X')$, where $X'$ is the base change $X \times_{\Spec R} \Spec R'$.
	This is because the projection $q\colon X' \to X$ is affine and thus $\O_X$-modules are the same as $q_*\O_X$-modules.
	A more modern way to spell this out would be to appeal to the Barr-Beck theorem, namely the underived and affine version of \cite[Theorem 4.7]{BFN}).
\end{rmk}
\subsection{Finiteness Conditions} 
\label{sub:finite}
We will be interested in characterising categorically quasi-coherent sheaves of finite type.
In an abstract abelian category $\cat{A}$ one can make the following definition.
We say that an object $M$ is of \emph{finite type} (or \emph{finitely generated}) if for any directed system of objects $(N_i)_{i \in I}$ the natural morphism
\begin{align*}
	\varinjlim \Hom\left( M, N_i\right) {\longto} \Hom(M, \varinjlim N_i).
\end{align*}
is injective.
We will check the equivalence with the standard definition in Proposition \ref{gerbe-ft}.
\begin{rmk}
	Although we shall not need it, we would like to point out a different characterisation of finite type quasi-coherent sheaves, which is meticulously explained in \cite{murfetab}.
	We might call an object $M \in \cat{A}$ categorically finitely generated if, for any family of subobjects $M_i \subset M$, such that $\sum_i M_i = M$ then there exists $i_0$ such that $M_{i_0} = M$.
	Using the fact that, on a quasi-compact and quasi-separated space, any quasi-coherent module is the direct limit of its finite type submodules \cite[Proposition 5.7.8]{RG}, one sees that categorically finitely generated objects coincide with finite type objects.
\end{rmk}

The other notion we need is that of \emph{compactness} for objects of the derived category.
An object $C$ of a triangulated category is said to be \emph{compact} if $\Hom(C,-)$ commutes with arbitrary coproducts.
For what we require we need to know the following facts hold for a quasi-compact and quasi-separated algebraic space: compact objects of the derived category are the same as perfect complexes \cite[\href{http://stacks.math.columbia.edu/tag/09m8}{Tag 09M8}]{stacks-project}; the (derived) dual of a perfect complex is perfect; the pullback of a perfect complex is perfect; the structure sheaf is a perfect complex.
	\section{Points and Graphs}
Let us fix a ground ring $R.$
We view an algebraic space $X$ over $R$ as a functor
\begin{align*}
	X\colon R\text{-}\Alg \longto \Set
\end{align*}
assigning to any $R$-algebra $A$ the set $X(A)$ of morphisms (over $R$) $\Spec A \to X$.
Thus we confuse $X$ with the moduli functor it represents, or in other words its \emph{functor of points}.
When $X$ is associated to the description of some objects, it is customary to call $X(A)$ the set of \emph{families} of objects parameterised by $\Spec A$.
The essence of this section is that, given an abelian category $\cat{C}$ over $R$, one can associate a moduli functor $\Pt_{\cat{C}}$.
Thus, if one wishes to define $\Pt_{\cat{C}}$ one needs to declare its values over all $R$-algebras, that is one needs to define what is a family of points.
To be precise, the functor $\Pt_\cat{C}$ takes values not in the category of sets but rather in the  $2$-category of groupoids, namely $\Pt_\cat{C}$ is a moduli \emph{stack}.

\subsection{Points} We now introduce our main definition, which mimics the properties enjoyed by graphs.
Let us fix a ground ring $R$.
Given an $R$-linear category $\cat{C}$ and a morphism $R \to A$ we shall denote the base change category $\cat{C}\otimes_R A$ also by $\cat{C}_A$.
If $P \in \cat{C}_{A}$ and $(-)\otimes_A P$ denotes the action of $\modo{A}$ on $P$, then a surjection $A \onto M$ is sent to an epimorphism $P \onto M \otimes_A P$.
This operation induces a well-defined function between equivalence classes of quotients of $A$ with equivalence classes of quotients of $P.$
\begin{center}
\begin{tikzcd}
	\left\{ \text{quotients of } A \right\} \arrow{r}{-\otimes_A P} & \left\{ \text{quotients of } P \right\}
\end{tikzcd}
\end{center}
\begin{defn}\label{definition}
	Let $\cat{C}$ be a cocomplete $R$-linear abelian category and let $A$ be an $R$-algebra.
	An \emph{$A$-point} of $\cat{C}$ (or a \emph{family} of \emph{pointlike objects} of $\cat{C}$ parameterised by $\Spec A$) is an object $P \in \cat{C}_A$ such that, for any morphism of $R$-algebras $A \to A'$ the following hold (denote by $P'$ the pullback of $P$ to $\cat{C}_{A'}$):
	\begin{enumerate}
		\item $P'$ is a finitely generated object of $\cat{C}_{A'}$;
		\item $P'$ is flat over $A'$;
		\item the functor
		\begin{align}\tag{$\circledcirc$}\label{compact}
			\RHom_{\cat{C}_{A'}}(-\ltensor A', P')\colon D(\cat{C})^{\opp} \to D(A')
		\end{align}
			sends compact objects to compact objects;
		\item the functor $\Phi = (-)\otimes_{A'} P'$ is fully faithful;
		\item $\Phi$ induces a bijection between equivalence classes of quotients of $A'$ and equivalence classes of quotients of $P'.$
	\end{enumerate}
	We form a groupoid $\Pt_\cat{C}(A)$ consisting of all $A$-points together with isomorphisms between them.
\end{defn}
From the definition we see immediately that the assignment $A \mapsto \Pt_\cat{C}(A)$ defines a prestack over $\Spec R$.
To obtain an honest set-valued functor, one relies on the usual trick.
Define
\begin{align}\label{P}
	\mathcal{P}_\cat{C} \colon R\text{-}\Alg &\longto \Set \\ \notag
						 A &\longmapsto \mathcal{P}_\cat{C}(A) = \left\{ P \in \Pt_\cat{C}(A) \right\}/_\sim
\end{align}
where $\sim$ stands for the equivalence relation which identifies two $A$-points $P_1$ and $P_2$ if there exists a line bundle $L$ on $\Spec A$ such that $L \otimes P_1 \simeq P_2$.
Notice that if $P$ is a pointlike object then $P \otimes L$ is also pointlike, as tensoring with $L$ is an autoequivalence of $\cat{C}_A$ and thus preserves pointlike objects.
As our definition of pointlike objects is purely in terms of the category theory of $\cat{C}$ the following lemma is merely an observation.
\begin{lem}\label{soup}
	Let $\cat{C}$ and $\cat{D}$ be two cocomplete $R$-linear abelian categories.
	If $\cat{C} \simeq \cat{D}$ as $R$-linear categories then $\Pt_\cat{C} \simeq \Pt_\cat{D}$ as prestacks over $\Spec R$ and $\mathcal{P}_\cat{C} \simeq \mathcal{P}_\cat{D}$ as presheaves over $\Spec R$.
\end{lem}

In some sense, the first four axioms defining points are technical, while the last one makes up the core.
Let us now unpack the definition in a geometric setting.
\begin{rmk}\label{pointsdowntoearth}
	Let $X$ be a quasi-compact and separated algebraic space over $R$.
	We want to understand the definition for $\cat{C} = \qc(X).$
	Given an $R$-algebra $A$ we have projection morphisms
	\begin{center}
	\begin{tikzcd}
	S \times X \arrow{d}{\pi}\arrow{r}{q} & X \\
	S &
	\end{tikzcd}$\quad\quad$
	where $S = \Spec A$ and the product $S \times X$ is taken over $\Spec R$.
	\end{center}
	An $S$-point of $\qc(X)$ is a quasi-coherent sheaf $P \in \qc(S \times X)$ such that the following properties hold (universally with respect to $S$):
	\begin{enumerate}
		\item $P$ is of finite type;
		\item $P$ is flat over $S$;
		\item for every compact complex $E \in D(X)$, the complex $\text{R}\pi_*\RlHom(\text{L}q^*(E),P)$ is compact;
		\item the natural map $M \to \pi_*\lHom_{S\times X}(P,P \otimes \pi^*M)$ is an isomorphism for all $M \in \qc(S)$;
		\item $\pi^*(-)\otimes P$ induces a bijection between closed subschemes of $S$ and quotients of $P$.
	\end{enumerate}
	A few comments are in order.
	The fact that $\qc(X) \otimes_R A \simeq \qc(S \times X)$ was already explained in Remark \ref{bchange}.
	The condition on compact complexes (see Subsection \ref{sub:finite}) is overkill and appears exclusively to ensure that $\pi_*P$ (corresponding to the case $E = \O_X$) is of finite type; see the step of Proposition \ref{succo} where we prove that $\pi_* P$ is a line bundle.
	By general category theory, $\Phi$ (which is a left adjoint) is fully faithful if and only the unit is an isomorphism.
	In practice we shall only use this for $M$ being the structure sheaf: $\O_S \simeq \pi_* \lHom(P,P)$.
	Notice again that if $P$ is an $S$-point, then so is $P \otimes \pi^*L$, for all line bundles $L$ on $S$.
\end{rmk}
We should mention that if one is solely interested in \emph{noetherian} algebraic spaces, then derived categories are unnecessary and one can replace \eqref{compact} by
$$\Hom(-\otimes A', P)\colon \cat{C} \to \modo{A'}$$
and compact objects of the derived category with finite type objects of the abelian category.
This will be made clearer in Remark \ref{noetherian}.
Our main goal now is to show that from $\Pt_{\qc(X)}$ one can recover $X$.
\subsection{Graphs} Let us now verify that graphs are indeed examples of pointlike objects.
Fix again a base ring $R,$ all fibre products are implicitly taken over $\Spec R.$
\begin{prop}\label{graphsarepoints}
	Let $X$ be a quasi-compact and separated algebraic space, $A$ be an $R$-algebra and $S = \Spec A$.
	Let $f\colon S \to X$ be a morphism.
	Consider the graph $\Gamma\colon S \to S \times X$.
	Then the structure sheaf of the graph $\Gamma_* \O_S$ is an $S$-point of the abelian category $\qc(X)$.
\end{prop}
The relevant diagrams are the following, notice that the square is cartesian.
\begin{center}
\begin{tikzcd}
S \times X \arrow{d}{\pi}\arrow{r}{q} & X \\
S &
\end{tikzcd}
$\quad\quad\quad$
\begin{tikzcd}
S \arrow{d}[swap]{f} \arrow{r}{\Gamma} & S \times X \arrow{d}{(f,\id)} \\
X \arrow{r}{\Delta} & X \times X
\end{tikzcd}
\end{center}
\begin{prf}
	Let $P$ be $\Gamma_* \O_S.$
	\begin{itemize}
		\item As $\Gamma$ is a closed immersion\footnote{It's worthwhile to point out that here we are using the separatedness of $X$.} we have that $P = \Gamma_* \O_S$ is of finite type.
		\item The functor $\pi^*(-)\otimes P$ $=$ $\pi^*(-) \otimes \Gamma_* \O_S$ $=$ $\Gamma_* \Gamma^* \pi^*$ $=$ $\Gamma_*$ is exact, therefore $P$ is flat over $S$.
		\item Let $E$ be a compact complex on $X$ and recall Subsection \ref{sub:finite}, which in particular says that in our setting compact objects are the same as perfect complexes \cite[\href{http://stacks.math.columbia.edu/tag/09m8}{Tag 09M8}]{stacks-project}.
		We have 
		$$\text{R}\pi_*\RlHom(\text{L}q^*E,P) = \text{R}(\pi \Gamma)_* \RlHom(\text{L}(q\Gamma)^*E,\O_S) = \RlHom(\text{L}f^*E,\O_S).$$
		As $E$ is perfect, the pullback $\text{L}f^*E$ is also perfect and so is its dual.
		\item As $\Gamma$ is a closed immersion $\Gamma_*$ is fully faithful.
	\end{itemize}
	The last axiom of being pointlike follows as we have already seen that the functor $\Phi = \pi^*(-)\otimes P$ is just $\Gamma_*,$ which satisfies the required property.
	Finally, to make sure these properties hold universally on $S$, notice that if $m\colon S' \to S$ is a morphism of affine schemes over $R$, then we have a cartesian diagram
	\begin{center}
	\begin{tikzcd}
		S' \arrow{r}{m} \arrow{d}[swap]{\Gamma'} & S \arrow{d}{\Gamma} \\
		S' \times X \arrow{r}{m_X} & S \times X \\
	\end{tikzcd}
	\end{center}
	where the vertical morphism on the left is given by the graph $\Gamma'$ of $fm$.
	As all morphisms in the diagram are affine, base change holds, and we have $m_X^*\Gamma_*(\O_S) = \Gamma'_* m^* \O_S = \Gamma'_* \O_{S'}.$
\end{prf}
\subsection{Points are Graphs} 
\label{sub:points_are_graphs}
The rest of this section is devoted to the crux of the paper: we prove that pointlike objects give rise to graphs.
\begin{prop}\label{succo}
	Let $X$ be an algebraic space over a base ring $R$ such that the structure morphism $X \to \Spec R$ is quasi-compact and separated.
	Let $A$ be an $R$-algebra, $S = \Spec A$ and $P$ an $S$-point of $\qc(X)$.
	Then there exists a unique morphism $f\colon S \to X$ of $R$-spaces such that, up to a twist by a unique line bundle on $S$, $P$ is the structure sheaf of its graph.
\end{prop}
We will divide the proof in small steps.
We start by fixing some notation.
Let $Z$ be the schematic support of $P$.
Denote its inclusion by $\iota$ and put $\rho = \pi \iota.$
\begin{center}
\begin{tikzcd}
	Z \arrow[right hook->]{r}{\iota} \arrow{dr}[swap]{\rho} & S \times X \arrow{r}{q} \arrow{d}{\pi} & X \\
	  &     S      &
\end{tikzcd}
\end{center}
\begin{rmk}\label{general.remarks}
	Before we start, a triple of general remarks.
	\begin{itemize}
		\item As $Z$ is the support of $P$, it follows $P \simeq \iota_* \iota^*P$ and, more generally, $E \otimes P \simeq \iota_* \iota^* E \otimes P$ for any $E \in \qc(S \times X)$.
		\item Additionally, as $P$ is flat over $S$ it follows that $\iota^*P$ is also flat over $S$.
		\item If $S' = \Spec A' \to S$ is a morphism and $P'$ is the pullback of $P$ to $S' \times X$ then, although \emph{a priori} the schematic support of $P'$ might not be equal to $Z \times_S S'$, the underlying topological spaces will be the same: $| \supp P' | = |Z \times_S S'|$ \cite[\href{http://stacks.math.columbia.edu/tag/07TZ}{Tag 07TZ}]{stacks-project}.
	\end{itemize}
\end{rmk}
Because it will come up as a key step below and because it makes for a good warmup, let us consider the case of a field.
\begin{lem}\label{field.case}
	Assume the ring $A$ to be a field.
	Then $\rho$ is an isomorphism.
\end{lem}
\begin{prf}
	As $A$ lacks any proper ideals, we deduce that any morphism $P \to Q$ is either the zero morphism or injective (this follows by using Axiom 5 and considering the image of $P$ in $Q$).
	As a consequence we have that if $P \onto Q$ is a surjection then either $Q = 0$ or $P \simeq Q.$
	Let $Z$ be the schematic support of $P$ and let us abuse notation by writing $P$ for the restriction of $P$ to $Z.$
	Let $I$ be a quasi-coherent ideal sheaf of $Z$ defining a closed algebraic subspace $W \subset Z$.
	Then $P/IP$ is either zero or isomorphic to $P.$
	We want to show that then $I = 0$ or $I = \O_Z$.
	We have
	\begin{align*}
		\left| \supp P/IP \right|  = \left| \supp \left( P \otimes \O_{Z}/I \right) \right| 
		= \left| \supp P \right| \cap \left| W \right| = \left| W \right|
	\end{align*}
	and therefore: if $P/IP = 0$, $W = \emptyset$ and $I = \O_Z$; if $P/IP = P$ then $I \subset \Ann(P) = 0$.

%
	
	Thus we see that $\O_Z$ has no non-trivial quasi-coherent ideal sheaves and thus is set-theoretically a singleton.
	We also know that there exists a dense open subset of $Z$ which is a scheme \cite[\href{http://stacks.math.columbia.edu/tag/06NH}{Tag 06NH}]{stacks-project}, hence $Z$ is a scheme.
	A scheme which has no non-trivial subschemes is the spectrum of a field.
	
	We now have a morphism of fields $A \to k,$ where $\Spec k = Z.$
	By Remark \ref{pointsdowntoearth} (4) the composition $A \to k \to \Hom_k(P,P)$ is an isomorphism.
	As $P$ has no non-trivial quotients we deduce that $P \simeq k$ and therefore $A \simeq \Hom_k(P,P) = k$.
	Thus, $\rho$ is an isomorphism.
\end{prf}
\noindent Let us now go back to the general case.

\smallskip
{\scshape The morphism $\rho$ is a universal homeomorphism.} 
Lemma \ref{field.case} and the third point of Remark \ref{general.remarks} imply that $\rho$ is universally bijective.
It is now enough to prove that $\rho$ is universally closed.
The functor $\Phi\colon \qc(S) \to \qc(S \times X)$ given by $\Phi(-) = \pi^*(-)\otimes P$ is the composition of $\rho^*$ followed by $\iota_*(-) \otimes P$.
Denote by $\cat{C}$ the essential image of $\Phi.$
By assumption $\Phi$ induces an equivalence between $\qc(S)$ and $\cat{C}.$
Moreover we know that $\Phi$ induces a bijection between the equivalence classes of quotients of $\O_S$ and of $P.$
If we combine this with the fact that $\supp P = Z$ we will see that $\rho^*$ induces (up to thickenings) a bijection between quotients of $\O_S$ and of $\O_Z.$
Drawing the following commutative diagram might be useful.
\begin{center}
\begin{tikzcd}
	\qc(S) \arrow{rr}{\Phi} \arrow{dr}[swap]{\rho^*} & & \cat{C} \\
	 	&	\qc(Z) \arrow{ur}[swap]{\iota_*(-)\otimes P}	&
\end{tikzcd}
\end{center}

As a linguistic matter, we say that two quotients $M \to M_1$, $M \to M_2$ of a fixed object $M$ are isomorphic \emph{as quotients} if they are isomorphic as objects in the under-category of $M$ (or, in other words, if they have the same kernel).
\begin{lem}
	Let $Y$ be an algebraic space, $F \in \qc(Y)$ a quasi-coherent sheaf of finite type such that $\supp F = Y$, i.e.~the schematic support of $F$ is equal to $Y$.
	Let $I,J$ be two quasi-coherent ideal sheaves defining two subspaces $Y_I, Y_J \subset Y.$
	If $F/IF \simeq F/JF$ as quotients of $F$ then the underlying topological spaces $\left| Y_I \right| = \left| Y_J \right|$ are the same.
	In other words, the operation $-\otimes F$ induces an injection (up to thickenings) from quotients of $\O_Y$ to quotients of $F$.
\end{lem}
\begin{prf}
	Notice that $F/IF \simeq F/JF$ as objects under $F$ if and only if $IF = JF.$
	We have
	\begin{align*}
		\left| Y_I \right| = \left| \supp F \right| \cap \left| Y_I \right| = \left| \supp F/IF \right|
		= \left| \supp F/JF \right| = \left| \supp F \right| \cap \left| Y_J \right|
	\end{align*}
	hence we are done.
	%
\end{prf}
Let us now show that the morphism $Z \to S$ is closed.
Let $\O_Z \onto \O_{Z'}$ be a quotient of $\O_Z$ defining a closed algebraic subspace $Z' \subset Z.$
Tensoring by $P$ gives a quotient of $P,$ which by assumption must lie in $\cat{C}$, the image of $\Phi.$
Thus there exists a quotient $\O_S \onto \O_{S'}$ such that $\pi^*(\O_{S'}) \otimes P = \iota_* (\O_{Z'} \otimes P)$, as quotients of $P$.
Call $\O_{Z''} = \rho^* \O_{S'}$ and notice that by construction $\iota_* (\O_{Z''} \otimes P) = \iota_* (\O_{Z'} \otimes P)$ as quotients of $P$.
By the previous lemma we thus have that the topological spaces $\left| Z' \right| = \left| Z'' \right|$ are the same.
Moreover, as $\O_{Z''} = \rho^* \O_{S'}$, we have $\left| Z' \right| = \rho^{-1}(\left| S' \right|)$.

As noted earlier, $\rho$ is universally bijective and, in particular, surjective.
Therefore one has than $\rho ( |Z'|) = \rho \rho^{-1}(|S'|) = |S'|$, thus implying the desired closedness of $\rho.$
%
To prove universal closedness of $\rho$, and not just closedness, we appeal once again to the third point of Remark \ref{general.remarks}.

\smallskip
{\scshape The sheaf $\pi_* P$ is a line bundle.} Here is where we use \eqref{compact}.
As $Z \to S$ is a universally closed and separated morphism with affine fibres it follows that it is affine \cite[Theorem 8.5]{rydh-approx}.\footnote{Many thanks are due to David Rydh for pointing out this key fact.}
As we have already observed that $\iota^*P$ is flat over $S,$ it now follows that the sheaf $\pi_*P  = \rho_* \iota^*P$ is flat over $S.$
To conclude that $\pi_* P$ is a bundle it suffices to show it is of finite presentation. 

Recall Subsection \ref{sub:finite}, in particular perfect and compact complexes coincide. In Axiom (3) we may take $E = \O_X$ to be the structure sheaf, which is perfect and hence compact.
This tells us that $\text{R}\pi_* P$ is compact and hence perfect.
Using the fact that $\rho$ is affine we have that $\text{R}\pi_* P = \text{R}(\pi \iota)_* \iota^* P = \rho_* \iota^* P = \pi_* P$ is a complex concentrated in degree zero.
As a perfect complex concentrated in degree zero is of finite presentation \cite[\href{http://stacks.math.columbia.edu/tag/066Q}{Tag 066Q}]{stacks-project} we obtain the claim.
Using Lemma \ref{field.case} we can compute the rank of $\pi_* P$ at each point of $S$, which is constantly one.
Hence, $\pi_*P$ is a line bundle.

\begin{rmk}\label{noetherian}
	If we assumed all our spaces to be noetherian one could slightly simplify the definition of a family of pointlike objects $P$. Instead of working with perfect complexes it would suffice to ask that for any finite type quasi-coherent sheaf $E$ on $X$, the module $\pi_*\lHom(q^*E,P)$ were of finite type.
	However, in the non-noetherian setting it is not clear whether graphs of morphisms would even satisfy this property!
	In fact, for a morphism $f\colon S \to X$ with graph $P = \Gamma_* \O_S$, we have $\pi_*\lHom(q^*E,P) = f^*E^\vee$.
	In general, the dual of a finite type module needn't be of finite type.
\end{rmk}

{\scshape The sheaf $P$ is (up to a twist) a graph.} 
Let us abuse notation and denote $\iota^*P$ by $P$.
As we know that $\rho_* P$ is a line bundle we have that the following sequence
\begin{align*}
	\O_S \to \rho_* \O_Z {\to} \End_{\rho_* \O_Z}(\rho_* P) \to \End_{\O_S}(\rho_* P) \simeq \O_S.
\end{align*}
composes to the identity, therefore the first morphism is injective and the third is surjective.
However, the third morphism is also injective (as it is a forgetful map) from which it follows that the second is surjective.
As the second is also injective (as $\supp P = Z$) it follows it is an isomorphism and therefore so is the first.
As $\rho$ is affine it follows that $Z \to S$ is an isomorphism.
If we denote by $f$ its inverse, with graph $\Gamma\colon S \to S \times X$, we have that $\Gamma_*\O_S = P \otimes \pi^* (\pi_* P)^\vee$.\hfill$\blacksquare$

%
%

\medskip
Adding together the results obtained so far allows us to prove Gabriel's theorem for algebraic spaces.
We present here the version suitable for audiences allergic to stacks and we then repeat the theorem in the next section using the language of gerbes, which is more natural from this paper's perspective.
\begin{thm}\label{maincomparison}
	Let $X$ be a quasi-compact and separated algebraic space over a ring $R$.
	Then, recalling the definition in \eqref{P}, the moduli functor $\mathcal{P}_{\qc(X)}$ is isomorphic to $X$.
	In particular, if $Y$ is another quasi-compact and separated algebraic space over $R$ then $\qc(X) \simeq \qc(Y)$ as $R$-linear abelian categories if and only if $X \simeq Y$ as algebraic spaces over $\Spec R$.
\end{thm}
\begin{prf}
	The proof has already been carried out in this section.
	What we've proved is that any $A$-point of $\qc(X)$ is, up to a twist of a line bundle on $\Spec A$ (see the last step of the previous proof) the graph of a morphism.
	That is, we've exhibited a functorial bijection between $\mathcal{P}_{\qc(X)}(A)$ and morphisms $\Spec A \to X$.
	Yet in other words, we've shown that the functors $\mathcal{P}_{\qc(X)}$ and $X$ are isomorphic.
	The second half of the theorem is immediate in light of Lemma \ref{soup}.
\end{prf}

	\section{Gerbes} 
\label{sec:the_main_result}
Let us fix yet once more a base ring $R$ and assume fibre products to be taken over $R$.
Given a cocomplete $R$-linear abelian category $\cat{C}$ we consider $\Pt_\cat{C}$ as a prestack over $\Spec R$, assigning to an $R$-algebra $A$ the groupoid $\Pt_\cat{C}(A)$ of $A$-points of $\cat{C}$.
We can soup up Lemma \ref{soup} as follows.

There is an action $\BGm \times \Pt_\cat{C} \to \Pt_\cat{C}$ given as follows.
Over a base $\Spec A$, $A$-linear autoequivalences act on $\Pt_{\cat{C}}(A)$.
In particular, tensoring with line bundles on $\Spec A$ realises an action of $\BGm$ on $\Pt_\cat{C}$.
%
%
\begin{lem}\label{mere}
	Let $\cat{C}$ and $\cat{D}$ be two cocomplete $R$-linear abelian categories.
	Then, if $\cat{C} \simeq \cat{D}$ as $R$-linear categories then $\Pt_\cat{C} \simeq \Pt_\cat{D}$ as prestacks over $\Spec R$.
	Moreover there is a canonical morphism
	\begin{align}\label{groupoid}
		\Equiv_R(\cat{C},\cat{D}) \longto \Iso_R^{\BGm}(\Pt_\cat{C},\Pt_\cat{D})
	\end{align}
	 from the groupoid of equivalences of $R$-linear categories $\Equiv_R(\cat{C},\cat{D})$ to the groupoid of isomorphisms of prestacks $\Iso^{\BGm}_R(\Pt_{\cat{C}},\Pt_{\cat{D}})$ compatible with the $\BGm$-action.
\end{lem}
The core of the previous section can be summarised as follows.
\begin{thm}\label{trivial}
	Let $X$ be a quasi-compact and separated algebraic space over the base ring $R$.
	The stack $\Pt_{\qc(X)}$ is isomorphic to $X \times \BGm.$
\end{thm}
\begin{prf}
	We define a map of stacks $X \times \BGm \to \Pt_{\qc(X)}.$ 
	If $A$ is an $R$-algebra and $S = \Spec A$, the groupoid $X \times \BGm(S)$ is given by the product of the set of morphisms $S \to X$, and the groupoid of line bundles $\BGm(S)$ on $S$.
	The natural map 
	$$(X \times \text{B}\G_\m)(S) \to \Pt_{\qc(X)}(S)$$
	sends the pair $(f\colon S \to X,L)$ to $(\Gamma_f)_*L$, which is equivalent to the structure sheaf of the graph of $f$, twisted by the line bundle $L$.
	Propositions \ref{graphsarepoints} and \ref{succo} show that a family of pointlike objects $P$ over $S$ is the same datum as a morphism $S \to \supp P \to S \times X \to X$ together with a line bundle over $S$, given by $\pi_* P$. Therefore, this map of groupoids is fully faithful and essentially surjective. This implies that we have defined an equivalence of stacks.
\end{prf}
As an immediate consequence we obtain again Theorem \ref{maincomparison} by killing off all automorphisms of $\Pt_{\qc(X)}$.
In other words we obtain the functor $\mathcal{P}_{\qc(X)}$, defined in \eqref{P}, by taking the quotient of $\Pt_{\qc(X)}$ by the action of the group-stack $\BGm$.

\subsection{General Gerbes} We remind the reader of the following equivalent concepts:
\begin{enumerate}
\item an element $\alpha \in H^2_{\text{\'et}}(X,\Gm)$,
\item a $\BGm$-torsor $\mathcal{X}$ on $X$,
\item a gerbe $\alpha\colon\mathcal{X} \to X$ with band $\BGm$ ($\Gm$-gerbe).
\end{enumerate}

We refer the reader to section 2.1 in \cite{lieb} for a definition of the concepts above.
Item (2) and (3) yield equivalent groupoids, while item (1) yields a set of objects, which is equivalent to the set of isomorphism classes of the groupoid of gerbes.
We say that a $\Gm$-gerbe $\alpha$ is \emph{neutral}, if $\alpha = 0$ in $H^2(X,\Gm)$.
A \emph{neutralization}, when it exists, is the choice of an equivalence $\mathcal{X} \cong X \times \BGm \cong [X/\Gm]$ as gerbes.
\begin{rmk}
	The group $H^2_{\text{\'et}}(X,\Gm)$ contains two distinguished subgroups: $\Br(X), \Br'(X)$, called respectively the \emph{Brauer group} and the \emph{cohomological Brauer group} of $X$.
	A \emph{Brauer class} is a class $\alpha \in H^2_{\et}(X,\Gm)$ which is induced by an Azumaya algebra, while $\Br'(X)$ is the torsion subgroup of $H^2_{\et}(X,\Gm)$.
	One has $\Br(X) \subset \Br'(X) \subset H^2_{\et}(X,\Gm)$ and in general both inclusions are strict.
\end{rmk}

Let $\alpha\colon \mathcal{X} \to X$ be a $\Gm$-gerbe on $X$, i.e.~\'etale locally on $X$ equivalent to $X \times \BGm$.
The abelian category $\qc(\BGm)$ decomposes as a product $\prod_{n \in \mathbb{Z}} qc(\BGm)_n$.
Each piece $\qc(\BGm)_n$ is equivalent to $\qc(\Spec R)$.
One obtains an analogous weight decomposition $$\qc(\mathcal{X}) \cong \prod_{n \in \mathbb{Z}} \QC(X,\alpha^n),$$ see Definition 2.1.2.1 in \cite{lieb}.
The category of $\alpha$-twisted sheaves $\qc(X,\alpha)$ is defined by means of this decomposition. A potential neutralization of $\alpha$ would induce an equivalence $\qc(X,\alpha) \cong \qc(X)$.
We can rephrase Theorem \ref{trivial} as saying that $\Pt_{\qc(X)}$ is the trivial $\Gm$-gerbe over $X.$
If one replaces $\qc(X)$ by $\qc(X,\alpha)$, the category of $\alpha$-twisted quasi-coherent sheaves (where $\alpha$ is a $\G_m$-gerbe), one still has a stack of points $\Pt_{\qc(X,\alpha)}$ and the theorem can be generalised.
Before we do so, we discuss the notion of twisted sheaves of \emph{finite type}.

In the following proposition we give a categorical characterization of finitely generated twisted sheaves.
We denote by $\lHom(-,-)$ the \emph{sheaf-hom} in the category $\O_X\text{-}\!\Mod$ of \'etale $\mathcal{O}_X$-modules on an algebraic space $X$.
For two $\alpha$-twisted sheaves $M$, $N$, the (untwisted) sheaf $\lHom(-,-)$ is well-defined. Characterizations of quasi-coherent sheaves of finite type as in the proposition below are well-known (see for instance \cite[Prop. 1.5.8]{bourbaki}). 

\begin{prop}\label{gerbe-ft}
Let $X$ be an algebraic stack over a ring $R$, which is quasi-compact and quasi-separated.
Let $\alpha$ be a $\Gm$-gerbe on $X$.
For an $\alpha$-twisted sheaf $M \in \qc(X,\alpha)$, the following assertions are equivalent:
\begin{enumerate}
\item $M$ is finitely generated,
\item for a directed system of $\alpha$-twisted sheaves $(N_i)_{i\in I}$ we have that $\varinjlim\lHom(M,N_i) \to \lHom(M,\varinjlim N_i)$ is injective in $\O_X\text{-}\!\Mod$,
\item for a directed system as above, we have that $\varinjlim\Hom(M,N_i) \to \Hom(M,\varinjlim N_i)$ is injective in $R\text{-}\!\Mod$.
\end{enumerate}
\end{prop}
\begin{prf}
	One should compare this proof with \cite[Prop. 71]{murfetmod}.
	
(1) $\Rightarrow$ (2): Since being finitely generated is an \'etale local notion, we may assume without loss of generality that $X$ is an affine scheme and that the gerbe $\alpha$ is neutral.
Thus $M$ is simply a module over an $R$-algebra $R'$, with generators $m_1,\dots,m_n \in M$. For an element $(f_i)_{i \in I} \in \varinjlim \Hom(M,N_i)$ to induce the zero map $M \to \varinjlim N_i$, we must have that for every $k \leq n$ there exists an $i_k \in I$ with $i \geq i_k$ implying $f_i(x_k) =0$.
Taking an upperbound $j \in I$ for the $i_k$, we have that for $i \geq j$ the map $f_i$ is zero.
Thus the system $(f_i)_{i \in I}$ is also zero as an element of $\varinjlim \Hom(M,N_i)$.

(2) $\Rightarrow$ (3): The global section functor $\Gamma(X,-)$ preserves directed colimits and maps injective maps of $\O_X$-modules to injective maps of abelian groups \cite[Tag 0739]{stacks-project}.

(3) $\Rightarrow$ (1): Let $f\colon U \to X$ be a quasi-compact and quasi-separated map of algebraic spaces.
The functor $i_*\colon \qc(U,i^*\alpha) \to \qc(X,\alpha)$ preserves directed colimits.
Let $(N_i)_{i\in I} \in \qc(U,i^*\alpha)$ be a directed system of $i^*\alpha$-twisted sheaves on $U$. Applying the adjunction between $i^*$ and $i_*$ to the injective map
$$\varinjlim \Hom_X(M, i_*N_i) \to \Hom_X(M,\varinjlim i_* N_i)$$
we obtain that 
$$\varinjlim \Hom_U(i^*M,N_i) \to \Hom_U(i^*M,\varinjlim N_i)$$
is injective.
If $U$ is an affine \'etale cover of $X$, neutralising the gerbe $\alpha,$ one reduces to the well-known case of rings \cite[Prop.~1.5.8]{bourbaki}.
%
%
\end{prf}
If one replaces ``finitely generated'' with ``finitely \emph{presented}'' in (1) and ``injective'' with ``bijective'' in (2) and (3), the theorem still holds with a slight modification to the proof.
We can now state the reconstruction theorem for categories of twisted sheaves.

\begin{thm}\label{twisted}
	Let $X$ be a quasi-compact and separated algebraic space over the base ring $R$ and let $\alpha\colon \mathcal{X} \to X$ be a $\Gm$-gerbe.
	The stack $\Pt_{\qc(X,\alpha)}$ is isomorphic to the gerbe corresponding to $\alpha.$
	
	Moreover, if $Y$ is a second quasi-compact and separated algebraic space equipped with a $\Gm$-gerbe $\beta\colon \mathcal{Y} \to Y$, \eqref{groupoid} translates into a morphism
	\begin{align*}
		\Equiv_R(\qc(X,\alpha),\qc(Y,\beta)) \longto \Iso^{\BGm}_R(\mathcal{X},\mathcal{Y})
	\end{align*}
	 from the groupoid of equivalences of $R$-linear categories $\Equiv_R(\qc(X,\alpha),\qc(Y,\beta))$ to the groupoid of $\BGm$-compatible isomorphisms of $R$-stacks $\Iso^{\BGm}_R(\mathcal{X},\mathcal{Y})$.
\end{thm}
\begin{prf}
	First off, let $X$ be an algebraic space and let $\alpha$ be a $\G_m$-gerbe given by $\alpha\colon \mathcal{X} \to X$.
	Denote by $\cat{T}$ the category $\qc(X,\alpha)$ of $\alpha$-twisted sheaves over $X$.
	Given $S \to \Spec R$ the base changed category $\cat{T}_S = \cat{T} \otimes_R S$ is isomorphic to $\qc(S \times X, q^*\alpha)$, where recall that $q\colon S \times X \to X$ is the second projection.
	
	We will define a morphism $\mathcal{X} \to \Pt_{\cat{T}}$. As above, let $S$ be an affine $R$-scheme. The groupoid of $S$-points $\mathcal{X}(S)$ is canonically equivalent to the groupoid of pairs $(f,\nu)$, where $f\colon S \to X$ is an $S$-point of $X$, and $\nu\colon X \times \BGm \cong f^*\mathcal{X}$ is a neutralization of the pullback of $\alpha$. The neutralization $\nu$ induces an equivalence of categories $\qc(S) \cong \qc(S,f^*\alpha)$, which we denote by $\mathbf{n}$.
	
	The graph of $f\colon S \to X$ gives rise to a pushforward functor $(\Gamma_f)_*\colon \qc(S,f^*\alpha) \to \qc(X \times S,\alpha)$. The map $\mathcal{X} \to \Pt_\cat{T}$ sends the $S$-point $(f,\nu)$ to the family of points $(\Gamma_f)_*\mathbf{n}(\mathcal{O}_S)$.
	Paralleling the proof of Proposition \ref{succo}, if $P \in \Pt_\cat{T}(S)$, we conclude that the schematic support $Z = \supp P$ of $P$ is universally homeomorphic to $S$ through the first projection $\pi\colon S \times X \to S$.
	This implies that $Z$ and $S$ have isomorphic \'etale sites \cite[Theorem 5.21]{rydh-sub}, therefore any open cover of $Z$ trivialising the gerbe structure on $P$ is the pullback of a cover on $S$.
	Passing to such a cover tells us that $Z$ and $S$ are isomorphic. Take now an open cover $U$ of $S$ which trivialises $f^*\alpha$.
	The untwisted case (Theorem \ref{trivial}) and faithfully flat descent imply now that the morphism $\mathcal{X} \to \Pt_\cat{T}$ is fully faithful and essentially surjective.
\end{prf}

A simple variant of the reasoning above allows us to generalise \cite[Theorem 6.1]{ant}, thereby removing the dependency on derived algebraic geometry from \cite{ant}.
\begin{cor}\label{underived}
Let $X$ be a quasi-compact and quasi-separated scheme, and $\alpha$, $\beta$ be two $\Gm$-gerbes on $X$. If $\qc(X,\alpha)$ and $\qc(X,\beta)$ are equivalent as Zariski stacks of $\mathcal{O}_X$-linear categories, then $\alpha \simeq \beta$.
\end{cor}

\begin{prf}
For every affine open $U \subset X$ we have an equivalence of $\Gamma(U,\mathcal{O})$-linear categories $\qc(U,\alpha|_U) \simeq \qc(U,\beta|_U)$. According to Theorem \ref{twisted}, this induces an equivalence of gerbes $\alpha|_U$, $\beta|_U$ on $U$. Since the intersection of finitely many affine open subsets is quasi-compact and separated, we obtain the descent data of an equivalence $\alpha \simeq \beta$.
\end{prf}


	\section{Autoequivalences}\label{autos}

It is folklore that for a smooth and projective variety $X$ the group of autoequivalences of the category $\coh(X)$ is isomorphic to the semidirect product $\Aut(X) \ltimes \Pic(X)$ (see for example the end of the proof of \cite[Corollary 5.24]{huybrechts}).
Here, $\Aut(X)$ denotes the group of automorphisms of $X$, and $\Pic(X)$ the group\footnote{For the categorically minded, both $\Aut(\qc(X))$ and $\Pic(X)$ are actually 2-groups.
As one usually does, however, we shall treat both as ordinary groups by taking isomorphism classes.
In the discussion which follows one can easily categorify our results, as the higher structure of $\Aut(\qc(X))$ and $\Pic(X)$ will always match up.} of line bundles on $X$.
We will generalise this result to quasi-compact and separated algebraic spaces \emph{flat} over a ring $R$.

For the last time, let us fix a ground ring $R$ and assume all fibre products to be taken over it.
If $X$ is an algebraic space over $R$ there is an obvious morphism
\begin{align}\label{picaut}
	\varsigma\colon\Aut_R(X) \ltimes \Pic(X) \to \Aut_R(\qc(X))
\end{align}
as automorphisms act on $\qc(X)$ via pullback and line bundles via tensor product.
The reason the product of groups above is semidirect is that the two actions do not commute: $f^*(M) \otimes L$ is in general different from $f^*(M \otimes L)$.
Going in the opposite direction, \eqref{groupoid} provides a morphism
\begin{align}\label{autpt}
	\Aut_R(\qc(X)) \longto \Aut_R^{\BGm}(\Pt_{\qc(X)}).
\end{align}
For $X$ quasi-compact and separated Theorem \ref{trivial} says that $\Pt_{\qc(X)} \simeq X \times \BGm$.
To connect this with our previous discussion we remind the reader of the following simple fact.
\begin{lem}
	If $X$ is an algebraic space over $R$ then $\Aut^{\BGm}_R(X \times \BGm) \simeq \Aut_R(X) \ltimes \Pic(X)$.
\end{lem}
\begin{prf}
	A $\BGm$-equivariant automorphism $X \times \BGm \to X \times \BGm$ is given by a morphism $$X \to X \times \BGm,$$ which induces an isomorphism $X \to X$. By the universal property of fibre products, it breaks therefore off into an isomorphism $X \to X$ and a map $X \to \BGm$, which corresponds to a line bundle on $X$.
\end{prf}
Thus the combination of \eqref{picaut}, \eqref{autpt} and the previous lemma tell us there is a right-split short exact sequence of groups
\begin{align*}
	1 \to Q_X \to \Aut_R(\qc(X)) \to \Aut_R(X) \ltimes \Pic(X) \to 1
\end{align*}
where $Q_X$ is just defined to be the kernel of \eqref{autpt}.
Assuming flatness of $X$ over $R$ one can show that $Q_X$ vanishes.
\begin{thm}\label{auto}
	Let $X$ and $Y$ be quasi-compact and separated algebraic spaces over a ring $R$ and assume either $X$ or $Y$ to be flat over a ring $R$.
	The natural morphism of Theorem \ref{twisted}
	$$\Equiv_R(\qc(X),\qc(Y)) \longto \Iso^{\BGm}_R(\Pt_{\qc(X)}, \Pt_{\qc(Y)})$$
	is an equivalence of groupoids.
	In particular, $\Aut_R(\qc(X)) \simeq \Aut_R(X) \ltimes \Pic(X).$
\end{thm}
The proof of this theorem is short, if one admits the existence of integral kernels for functors between derived $\infty$-categories.
Since these results require methods from derived algebraic geometry, the claim of elementariness in the abstract does not extend to this section.

\begin{prf}
	Assume $X$ to be flat over $R.$
	We denote by $\QC(X)$ the unbounded derived category of $\qc(X)$, which we view as an $R$-linear stable $\infty$-category.
	An equivalence $\psi\colon \qc(X) \to \qc(Y)$ extends automatically to an equivalence $\Psi\colon \QC(X) \to \QC(Y)$.
	We also denote by $\varphi\colon \Pt_{\qc(X)} \to \Pt_{\qc(Y)}$ the isomorphism corresponding to $\psi.$
	In \cite[Theorem 1.2 (2)]{BFN} (see also \cite{orlovk3}, \cite{toenmorita}) it is shown that any such $\Psi$ is given by a Fourier-Mukai transform, provided $X$ and $Y$ are \emph{perfect} in the sense of \cite[Definition 3.2]{BFN}.
	Notice that here we use the flatness assumption to ensure that the \emph{derived} fibre product of $X$ and $Y$ over $\Spec R$ is just the ordinary fibre product $X \times Y$.
	Any quasi-compact algebraic space with affine diagonal is perfect thanks to \cite[Proposition 3.9]{BFN}, \cite[\href{http://stacks.math.columbia.edu/tag/09M8}{Tag 09M8}]{stacks-project} and \cite[\href{http://stacks.math.columbia.edu/tag/09IY}{Tag 09IY}]{stacks-project}.
	Let $K \in \QC(X \times Y)$ be the integral kernel representing $\Psi$.
	Then $\O_X \boxtimes K \in \QC(X \times X \times Y)$ is the kernel of a Fourier-Mukai equivalence $\Psi_X$ between $X \times X$ and $X \times Y$, seen as spaces over $X$.
	A standard computation with Fourier-Mukai functors shows that $K = \Psi_X(\Delta_* \O_X)$, where $\Delta_* \O_X$ is the image of the structure sheaf of the diagonal of $X \times X$.
	Since $\Psi$ comes from the equivalence $\psi$, we see that $K = \varphi(\Delta_*\O_X) \in \Pt_{\qc(Y)}(X)$ is an $X$-family of points of $\qc(Y)$.
	Accordingly, we see that $\Psi$ (and thereby $\psi$) is completely determined by $\varphi.$
\end{prf}

\begin{thm}\label{tauto}
Let $X$ and $Y$ be quasi-compact and separated algebraic spaces over $R$, $\alpha\colon \mathcal{X} \to X$, $\beta\colon \mathcal{Y} \to Y$ two $\G_m$-gerbes on $X$, respectively $Y$.
Assume either $X$ or $Y$ to be flat over $R$.
Then, the natural map 
$$\Equiv_R(\qc(X,\alpha),\qc(Y,\beta)) \longto \Iso_R^{\BGm}(\mathcal{X},\mathcal{Y})$$
is an equivalence.
\end{thm}
\begin{prf}
	Perfectness of $\mathcal{X}$ boils down to knowing whether the unbounded derived category $\QC(X,\alpha)$ is compactly generated.
	In forthcoming work of Hall-Rydh this is shown to be true for a broad class of spaces, including $\Gm$-gerbes over tame Artin stacks.
	Hence we can apply the formalism of \cite[Thm 1.2(2)]{BFN}.
To a colimit-preserving functor $\QC(X,\alpha) \to \QC(Y,\beta)$ we associate the colimit-preserving functor
\begin{equation}\label{silly}
\QC(\mathcal{X}) \cong \prod_{n \in \mathbb{Z}} \QC(X,\alpha^n) \to \prod_{n\in \mathbb{Z}} \QC(Y,\beta^n) \cong \QC(\mathcal{Y}),
\end{equation}
given by the zero functor for $n \neq 1$.
By virtue of \cite[Thm 1.2(2)]{BFN} there is a complex of sheaves $K \in \QC(X \times Y,\alpha^{-1} \boxtimes \beta)$ which gives rise to the functor of \eqref{silly}. The identity functor from $\QC(X,\alpha)$ to $\QC(X,\alpha)$ is given by endowing the sheaf $\Delta_*\mathcal{O}_X$ with the structure of an $\alpha^{-1}\boxtimes \alpha$-twisted sheaf.
This is possible, since the gerbe $\Delta^*(\alpha^{-1}\boxtimes\alpha)$ is neutral.
One can now repeat the proof of Theorem \ref{auto}.
\end{prf}

	\section{Remarks} 
\label{sec:comments}
\begin{rmk}
	In \cite{gaitsgory} it is shown that there is also a functioning theory of descent for abelian categories.
	So far we have worked with a base ring $R.$
	If we replace $\Spec R$ by a general scheme or (higher) stack $S$ one can speak about $S$-linear abelian categories and sheaves (or rather $2$-sheaves) of abelian categories over $S$.
	The prototypical example being again $\qc(X)$ for a stack $X$ over $S$.
	It is not surprising that in the affine case $S = \Spec R$, a sheaf of abelian categories is equivalent to its global sections category.
	The following weak result can then be bootstrapped from our main result.
	\begin{cor}\label{morita}
		If $S$ is an algebraic stack and $X$ and $Y$ are stacks over $S,$ which are relatively representable by quasi-compact and separated algebraic spaces, then $X \simeq Y$ as $S$-stacks if and only if $\qc(X) \simeq \qc(Y)$ as \emph{sheaves} of $\O_S$-linear abelian categories.
	\end{cor}
\end{rmk}

	In a similar vein, one can study sheaves of abelian categories over a stack $X$ which are locally the category of quasi-coherent sheaves.
	This produces a \emph{Morita theory} for sheaves of abelian categories.
	Simply put, this is a higher version of the statement that a line bundle is the same thing a $\Gm$-torsor.
	The way we prove this is by realising that the cocycle description for either is identical: an automorphism of a line bundle is given by an element of $\O^\times$.
	
	Let $X$ be a quasi-compact and quasi-separated algebraic stack over a ring $R$.
	To be safe, we note that if $U \to X$ and $V \to X$ are quasi-compact and quasi-separated algebraic stacks mapping to $X$ then so is $U \times_X V$.
	A $\BGm$-torsor, or a $\Gm$-gerbe, can be given a cocycle description as follows.
	If $\{U_i \to X\}_i$ is a trivialising smooth cover (which we can assume to be affine) we are left with the datum of: $U_i \times \BGm$ on each $U_i$, a $\BGm$-equivariant equivalence over each double intersection $U_{ij}$ (i.e.~a line bundle $L_{ij}$ on $U_{ij}$) and 2-isomorphisms between the composition of the $1$-isomorphisms above over triple intersections $U_{ijk}$ (equivalently, an isomorphism $L_{ij}| \otimes L_{jk}| \simeq L_{ik}|$, where $|$ stands for ``restriction to $U_{ijk}$'') satisfying the quadruple cocycle identity in $\O^\times_{U_{ijkl}}$.
		
		On the other hand, one can apply Theorem \ref{twisted} to obtain the very same description for sheaves of abelian categories.
		Let $\mathcal{C}$ be a sheaf of $\O_X$-linear abelian categories over $X$ for the smooth topology.
		We say that $\mathcal{C}$ is \emph{invertible}\footnote{This should be compared with \cite[Defn 2.1.4.1 and Prop 2.1.5.6]{max}.} if there exists a smooth cover $\{U_i \to X\}_i$ together with $\O_{U_i}$-linear equivalences $\mathcal{C}(U_i) \simeq \qc(U_i)$.
		By Theorem \ref{auto}, $\Aut_{\O_{U_{ij}}}(\qc(U_{ij})) \simeq \Pic(U_{ij})$ (and similarly for $n$-tuple intersections),\footnote{As we are working over $U_{ij}$, $\Aut_{U_{ij}}(U_{ij})$ is trivial and the flatness is automatically satisfied.} so that the cocycle description of $\mathcal{C}$ is identical to that of a $\Gm$-gerbe over $X$.
		In other words we have the following result.
		\begin{cor}
			Let $X$ be a quasi-compact and quasi-separated algebraic stack.
			The $2$-groupoid\footnote{I.e.~a $2$-category where morphisms and $2$-morphisms are invertible.} of $\BGm$-torsors on $X$ is equivalent to the $2$-groupoid of invertible sheaves of $\O_X$-linear abelian categories on $X$.
		\end{cor}

\begin{rmk}
	Another class of spaces where one might consider a generalisation of Gabriel's theorem is that of \emph{derived schemes}.
	In this context one must consider the whole triangulated (or rather stable $\infty$-)category $\QC(X)$ of unbounded quasi-coherent modules, as it does not arise as a derived category.
	It is well known that $\QC(X)$ as a triangulated category does not recover $X$, for example when $X$ is an abelian surface \cite[Theorem 2.2]{mukai}.

	However, $\QC(X)$ comes with a t-structure, whose heart is $\qc(X_0)$, the category of quasi-coherent sheaves on the underlying underived scheme.
	One might then consider the pair $(\QC(X),\qc(X_0))$, consisting of a triangulated category equipped with a t-structure.
	Unfortunately, even in this sense, the reconstruction theorem fails already for derived affine schemes.
	In fact, there are commutative differential graded algebras which are quasi-isomorphic as differential graded algebras but not as \emph{commutative} differential graded algebras.\footnote{We would like to thank Bertrand Toën for pointing this out.}
\end{rmk}

	\bibliographystyle{fabio.bst}
	\bibliography{biblio.bib}
\end{document}